\documentclass{article}

\begin{document}

\title{Germain and Her Fearless Attempt to Prove Fermat's Last Theorem}
\author{Dora Musielak}
\date{}

\maketitle

\begin{abstract}
Two centuries ago, Sophie Germain began to work on her “grand plan to prove the theorem of Fermat,”
 the famous conjecture that $x^n+y^n=z^n$ is impossible for nonzero integral values of $x$, $y$, and $z$, when $n>2$.  
At that time, this was an open question since nobody knew whether Fermat's assertion was true. 
Euler had proved it for $ n=3$ and $n=4$. However, no one else had demonstrated the general case.
Then Sophie Germain valiantly entered the world of mathematics in 1804, reaching out to Gauss
 (writing under the assumed name Monsieur Le Blanc) boldly stating that she could do it. 
Eventually, Germain conceived a formidable plan for proving Fermat's Last Theorem in its entirety, 
and in the process she obtained proofs of Case 1 for particular families of exponents. 
Her efforts resulted in Sophie Germain’s Theorem that proves Case 1 of FLT for an odd prime exponent $p$ whenever $2p+1$ is prime. 
Today, a prime $p$ is called a Sophie Germain prime if $2p+1$ is also prime. 
It remains an unanswered question whether there are an infinite number of 
Sophie Germain primes. But there is more that Germain did in number theory,
much of which was veiled by the mathematicians with whom she shared her work. 
This article provides historical details of Germain's efforts, written with the sole intention 
of paying homage to the only woman mathematician who contributed to proving the most
famous assertion of Fermat.
\end{abstract}

$\hskip1.0in$              {\it ``Time keeps only the works that defend against it.''} \

$\hskip3.15in$             {Sophie Germain}


\section{Introduction}

For two centuries, historians of mathematics overlooked the contributions of Sophie Germain to the proof of Fermat's Last Theorem (FLT). In 1825, Adrien-Marie Legendre relegated her effort to a footnote in the Second Supplement of his {\it Essai sur la th\'eorie des nombres}. Yet, Germain had done much more than the proof that Legendre credited her. And he knew that very well because she communicated her results to him through letters. Gauss knew that as well, for in a letter in 1819, Germain sent him a treasure trove outlining her plan for the general proof of the {\it c\'el\`ebre \'equation de Fermat}, as she called it.

In order to better appreciate Germain's attempt to prove FLT, we first need to present a perspective on the background that sustained her work. Let's start with Fermat, who is responsible for the puzzling assertion that became so famous.{\footnote{It is known as Fermat's Last Theorem because it was the last of his assertions that remained unproved.}}

Pierre de Fermat is considered the father of the modern theory of numbers. A lawyer, Fermat was councillor for the parliament of Toulouse in 1631. Fermat devoted his leisure time to mathematics, making important contributions along the way. And although he communicated his results with other scholars, Fermat was secretive about his analytical methods or proofs. Fermat owned a copy of Bachet's {\it  Diophantus}{\footnote{In 1612, Bachet de M\'eziriac (1581-1638) published {\it Probl\'emes plaisants et d\'electables qui se font par les nombres}, and in 1621 a Greek edition of {\it Diophantus} with notes.}},  in which he wrote numerous marginal notes to record his ideas. 

Five years after his death, his son published a new edition of {\it Diophantus}, incorporating Fermat's notes as an appendix. In the margin of one page, Fermat had written, ``It is impossible to divide a cube into two cubes, or a fourth power into two fourth powers, or in general any power beyond the square into powers of the same degree; of this I have discovered a very wonderful demonstration {[\it demonstrationem mirabilem sane detexi}]. This margin is too narrow to contain it.''{\footnote{Fermat did leave a proof for $n = 3$. However, it is unlikely that he had a valid proof for the general case. Mathematics historians have written that FLT dates from about the year 1637. However, Fermat died in 1665. Why did Fermat not communicate the proof to anyone in all those many years since he conjectured it?}  What Fermat meant is that there is no positive integer $n>2$ for which 

\begin{equation}
    \label{simple_equation}
     x^n + y^n = z^n
\end{equation}
can be solved in terms of positive integers $x$, $y$, $z$.

After Fermat, Leonhard Euler gave luster and depth to number theory. 
It was Christian Goldbach who (in a letter dated December 1729) introduced 
the twenty-two-year old Euler to the work of Fermat.{\footnote{Other theorems due to
Fermat on numbers were published in his {\it Opera varia} and in Wallis's {\it Commercium epistolicum} of 1658.}}  In the postscript of that letter,{\footnote{Euler correspondence with Goldbach in Fuss, P.-H., ed., {\it Correspondance Math\'ematique et Physique}, Imp. Acad. Sci., St. Petersburg, 1843, vol. I, {\it see} [F], Letter II, p. 10.}  Goldbach wrote:

P.S. {\it Notane Tibi est Fermatii observatio omnes numeros hujus formulae} $2^{2^{x-1}}+1$ {\it nempe} 3, 5, 17, {\it etc.
esse primos, quam tamen ipse fatebatur sc demonstrare non posse, et post eum nemo, quod sciam, demonstravit.}
 [P.S. Note Fermat's observation that all numbers of this form $2^{2^{x-1}}+1$, that is 3, 5, 17, etc. are primes, but he 
 confessed that he did not have a proof, and after him, no one, to my knowledge, has proved it.]
 
The terse postscript sparked the young man's interest. Euler responded that he could not find a proof either, [{\it Nihil prorsus invenire potui, quod ad Fermatianam observationem spectaret.}]{\footnote{{\it Nihil prorsus invenire potui, quod ad Fermatianam observationem spectaret. Sed nondum prorsus persuasus sum, quomodo sola inductione id inferre legitime potuerit, cum certus sim ipsum numeris in formula $2^{2^x }$ loco $x$ substituendis nee ad senarium quidem pervenisse.} Fuss, p. 18, {\it see} [F].}}  However, two years later, he disproved Fermat's assertion. In a 1732 paper titled, ``Observations on a theorem of Fermat and others on looking at prime numbers,''{\footnote{Euler, L., {\it Observationes de theoremate quodam Fermatiano aliisque ad numeros primos spectantibus}, Commentarii academiae scientiarum Petropolitanae 6, 1738, pp. 103-107. Paper E26, available at The Euler Archive: http://eulerarchive.maa.org/.}}  Euler stated that  ``{\it Est enim} $2^{2^5} +1=2^{32}+1=4294967297$,'' a number which he determined is not prime since it is divisible by 641. 
And so, Euler began proving and disproving Fermat's assertions. Then he tried the proof of (1).

In a memoir{\footnote{Euler, L., {\it Theorematum quorundam arithmeticorum demonstrationes} (The proofs of some arithmetic theorems), presented to the St. Petersburg Academy on June 23, and August 16 (Additions), 1738. E98 in The Euler Archive.}}  published in 1738, Euler used the method of infinite descent to prove that  $x^4+y^4 \neq z^4$. This method, known since antiquity, is described as follows: suppose that a solution of a given problem is possible in positive integers: then one can show how to derive from it a solution in smaller positive integers, and so on. But since this process cannot go on indefinitely, we reach a contradiction and thus show that no solution is possible. Euler also gave a simpler proof for $n=4$ in his book {\it Algebra},{\footnote{Euler, L., {\it Vollstandige Anleitung zur Algebra} (Complete instruction in algebra), {\it Opera Omnia}: Series 1, Volume 1. Papers E387, E388 in The Euler Archive. }}  published in 1770, where we also find his proof for $n=3$. 

In 1798, a year before the end of the French Revolution, Adrien-Marie Legendre published {\it Essai sur la th\'eorie des nombres.}{\footnote{Legendre presented his {\it Essai} at the meeting of the Class of Mathematics on 13 July 1798. {\it S\'eance du 25 Messidor An} 6, {\it see} [IP], Tome I.}}  In this treatise, Legendre compiled his researches on number theory and added the numerous but scattered fragments on the theory of numbers due to his predecessors (Diophantus, Vi\`ete, Bachet, Fermat, Euler, Lagrange), arranging the body of work into a systematic whole. The {\it Essai} was preceded by another memoir he published in 1785 with the title ``Research of Indefinite Analysis,''{\footnote{Legendre, A.-M., {\it Recherches d'analyse ind\'etermin\'ee}, Paris, 1785.}}  which refers mainly to the properties of numbers and where we find Legendre's incomplete demonstration of the law of quadratic reciprocity issued by Euler. 

Three years later, Johann Carl Friedrich Gauss published his {\it Disquisitiones Arithmeticae},{\footnote {Gauss, C. F., {\it Disquisitiones arithmeticae}. Brunswick, July 1801.}}  a treatise that revolutionized the study of number theory. As Gauss wrote in the preface, the science of Arithmetic was much indebted to Fermat, Euler, Lagrange and Legendre who had ``opened the door to this divine science and discovered a mine of inexhaustible riches in it.'' Gauss stated the scope as follows: ``The inquiries which this volume will investigate pertain to that part of Mathematics which concerns itself with integers.'' Written in Latin, Gauss addressed in the {\it Disquisitiones} both elementary number theory and parts of the area of mathematics now called algebraic number theory. In addition to reconciling results in number theory obtained by the mathematicians he praised, Gauss added many profound and original results of his own. 

The treatises of Legendre and Gauss appeared in print right at the moment when Sophie Germain crossed the threshold into mathematics. Through their work, both number theorists gave her the tools and technical background, and also the inspiration to tackle one of the most intriguing and challenging problems in mathematics. I must point out here that, when Euler wrote to Goldbach in 1753 to share that he had proved FLT in the cases $n=3$ and $n=4$, he stated that ``for these two cases are so different from each other, that I see no possibility of deriving from them a general proof for $a^n+b^n \neq c^n$ if $n>2$.''{\footnote{Fuss's {\it Correspondence}, Letter CLV, p. 618, {\it see} [F].}} Euler added: ``But one can see quite clearly by extension that the greater the $n$, the more difficult the demonstrations must be.''{\footnote{Ibid.}  Euler was right.

That is how the theorem of Fermat stood when Sophie Germain encountered it. A proof of the general case of FLT would require development of new mathematical methods, and more than two centuries of efforts to be discovered. But Sophie Germain did not know that, so she attempted it, perhaps driven to solve the puzzle by the sheer simplicity of the equation (1) and lured by Fermat's assertion that he had discovered ``a very wonderful demonstration.''

\section{Who was Sophie Germain?}

Sophie Germain was born in 1776 in Paris. Although there are no historical records to verify it, I believe Germain taught herself elementary mathematics in her teenage years.{\footnote{Musielak, D., {\it Sophie's Diary}. Mathematical Association of America, MAA Spectrum Book Series, ISBN 1-4184-0812-3, Second Edition, 2012.}}  Then, at age twenty-two, Sophie Germain reached out to Joseph-Louis Lagrange, hiding behind a man's name (Monsieur Le Blanc), to submit her analysis, allegedly prompted by his lecture notes at the \'Ecole Polytechnique. 

Soon, Lagrange discovered her true identity and this caused a stir in Paris. While some members of the Parisian Academy of Sciences sent books and offered to help Germain in her studies, others were fascinated by the notion of a lady mathematician and pestered her with silly public adulation that did not sit well with her. Astronomer Lalande, for example, visited Germain at her home and said something that insulted her intelligence (he later apologized). A lyricist composed and published a poem that caused Germain and her family much embarrassment.{\footnote{Stupuy, H., {\it Oeuvres philosophiques de Sophie Germain, suivies de pens\'ees et de lettres in\'edites et pr\'ec\'ed\'ees d'une \'etude sur sa vie et ses oeuvres}, Paris, 1879, {\it see} [SH].}} Such experiences led Sophie Germain to believe that scholar women were ridiculed rather than taken seriously.

There is no historical record to help us state with certainty that it was actually Legendre who came to her aid, won her trust, and became her mentor. What is evident is that the humble, well-mannered mathematician did take Sophie Germain seriously. I believe Germain used Legendre's work to learn number theory and introduce her to Fermat's theorems. In fact, the year when Legendre published the {\it Essai} (1798) coincides with the time when Germain was discovered as a clandestine student of mathematics.{\footnote{Musielak, D.,{\it Sophie Germain}, 2020, Springer Biography, {\it see} [MD].}}  

Through correspondence, which she initiated, Legendre became Germain's first teacher and advocate. Number theory was reemerging as an important branch of mathematics, thanks in part to the work of Legendre, and especially the {\it Disquisitiones Arithmeticae} of Gauss. 

On 26 January 1802, Legendre communicated to the French Academy the discovery of Gauss's work. In all likelihood, Legendre shared this news with Germain and probably encouraged her to study it. By 1804, she was intent on proving the famous assertion of Fermat. She was twenty-eight years old and considered herself an enthusiastic amateur ({\it amateur enthousiaste}) mathematician. Whether she was persuaded to read it or she discovered the {\it Disquisitiones} on her own, Germain wrote to Gauss saying that his work had been ``the object of her admiration and her studies'' for a while.{\footnote{Stupuy, H., {\it Oeuvres philosophiques de Sophie Germain, suivies de pens\'ees et de lettres in\'edites et pr\'ec\'ed\'ees d'une \'etude sur sa vie et ses oeuvres}, Paris, 1879, Letter VI, pp. 254-258, {\it see} [SH].}}  At that time, Fermat's assertion, still not known as Fermat's Last Theorem, was not the bewildering conjecture that became legendary in the twentieth century. Yet, the theorem caught Germain's attention and proving it became her obsession. 

Nevertheless, still hurt by her experience with condescending male intellectuals, Germain resorted to pretending she was a man when writing to Gauss. Eventually, he discovered that Monsieur Le Blanc was actually Mademoiselle Germain.{\footnote{Musielak, D., {\it Sophie Germain}, chapter 3, {\it see} [MD].}} He was rather pleased and, from then on, he came to consider Sophie as a friend.

\section{Germain's First Attempt to Prove FLT}

In her first letter to Gauss, dated 21 November 1804, Sophie Germain stated that she could prove the th\'eor\`eme de Fermat, that $x^n+y^n=z^n$ is impossible if $n=p-1$, where $p$ is a prime of the form $8k+7$. Her statement suggests that Germain had acquired a good understanding of number theory, enough to perceive the full depth of the unproven assertion. To prove Fermat's conjecture in general, it is sufficient to demonstrate the impossibility of

\begin{equation}
    \label{simple_equation}
     x^p + y^p = z^p
\end{equation}
in non-zero integers $x$, $y$, and $z$, for any odd prime $p>3$. In fact, Legendre proved in 1823 the case for $p=5$.

The first paragraphs of that important first letter to Gauss represent Germain's eager attempts to be recognized by the German mathematician who was at that moment ``placed in the ranks of the best geometers.''{\footnote{``Geometer'' is equivalent to ``pure mathematician.'' In a letter dated 31 May 1804, Lagrange wrote to Gauss: ``Your Disquisitiones have placed you immediately in the ranks of the best geometers.''}}  Note her statement in the second paragraph concerning ``the famous equation of Fermat'':\\

{\it Monsieur},

{\it Vos Disquisitiones arithmeticae sont depuis longtemps l'objet de mon admiration et mes \'etudes. Le dernier chapitre de ce livre renferme, entre autres choses remarquables, le beau th\'eor\`eme contenu dans l'\'equation $4 ({x^n-1})/(x-1)=Y^2 \pm nZ^2$; je crois qu'il peut \^etre g\'en\'eralis\'e ainsi $4 ((x^{n^s}-1))/(x-1)=Y^2 \pm nZ^2$, $n$ \'etant toujours un nombre premier et $s$ un nombre quelconque. Je joins \`a ma lettre deux d\'emonstrations de cette g\'en\'eralisation. Apr\`es avoir trouv\'e la premi\`ere j'ai cherch\'e comment la m\'ethode que vous avez employ\'ee art. 357 pouvait \^etre appliqu\'ee au cas que j'avis à considérer. J'ai fait ce travail avec d'autant plus de plaisir qu'il m'a fourni l'occasion de me familiariser avec cette m\'ethode qui, je n'en doute pas, sera encore dans vous mains l'instrument de nouvelles d\'ecouvertes. J'ai ajout\'e \`a cet article quelques autres consid\'erations. 

La derni\`ere est relative \`a la c\'el\`ebre \'equation de Fermat $x^n+y^n=z^n$ dont l'impossibilit\'e en nombres entiers n'a encore \'et\'e d\'emontr\'ee que pour $n=3$ et $n=4$. Je crois \^etre parvenu \`a prouver cette impossibilit\'e pour $n=p-1$, $p$ \'etant un nombre premier de la forme de $8k+7$. Je prends la libert\'e de soumettre ces essais \`a votre jugement, persuad\'e que vous ne d\'edaignerez pas d'\'eclairer de vos avis un amateur enthousiaste de la science que vous cultivez avec de si brillants succ\`es.}\\

Six months later, Gauss responded, but not with the feedback Germain expected. After apologizing for the late reply and reaffirming his delight for the {\it  recherches arithm\'etiques} to which he had devoted the most beautiful part of his youth ({\it la plus belle partie de ma jeunesse}), Gauss explained that his work was now devoted to astronomy, research which absorbed all his time. Nonetheless, thinking that his correspondent was Monsieur Le Blanc, Gauss addressed another enclosed proof regarding prime numbers for which 2 is a certain power residue or nonresidue, saying that it was weak.{\footnote{{\it Surtout votre nouvelle d\'emonstration pour les nombres premiers, dont 2 est r\'esidu ou non r\'esidu [quarr\'e], m'a extr\^emement plu ; elle est tr\`es fine, quoiqu'elle semble \^etre isol\'ee et ne pouvoir s'appliquer \`a d'autres nombres.} [``I liked especially your new proof regarding prime numbers for which 2 is a [certain power] residue or non residue; it is fine but it seems to be an isolated case not applicable to other numbers.'']. Letter VII dated 16 June 1806 in Stupuy's Correspondence, pp. 258-261, {\it see} [SH].}} The theory of power residues, that Gauss championed in 1801, would play a central role in Germain's approach towards a proof of FLT.

Gauss did not elaborate on her proof related to the impossibility of Fermat's theorem with the prime exponent $n=p-1$, which was the most important to Germain and also to the advancement of the analytical procedures in number theory at that time.{\footnote{There is no historical record to suggest that Lagrange, Legendre or Gauss ever made attempts to prove the general case of FLT (for all $n$).}} 

Germain wrote to Gauss on 20 February 1807, stating that if the sum of the $n$th powers of any two numbers is of the form  $h^2+nf^2$, the sum of these two numbers is of that form. Promptly, Gauss replied{\footnote{Gauss letter is dated 30 April 1807. Published by Stupuy in 1879, {\it see} [SH], p. 276.}} that this is false, and he added an example where the rule fails.

Sophie Germain quickly adopted the congruence notation that Gauss introduced in 1801, and she used it to articulate her mathematical ideas. Germain continued sending samples of her work to Gauss, asking for his critique. And even though number theory was his first love ({\it mon arithm\'etique ch\'erie}),{\footnote{``... {\it qui ont \'et\'e la source de mes jouissances les plus d\'elicieuses et qui me seront toujours plus ch\'eres qu’aucune autre science}.'' Gauss letter to Germain dated 16 June 1805, {\it see} [SH].}}  Gauss was not keen to engage on her attempts to prove the famous theorem of Fermat. 

On 12 January 1807, Delislc Poulet presented to the French Academy the translation of Gauss's {\it Disquisitions} into French.{\footnote{Poulet was a professor at the Lyc\'ee d'Orl\'eans. The translated {\it Disquisitions} became {\it Recherches arithm\'etiques}.}}  Shortly after, in a letter dated 20 February, Germain confessed to Gauss that Monsieur Le Blanc was actually her.{\footnote{Musielak, D., {\it Sophie Germain}, chapter 3, {\it see} [MD].}} He was delighted. And although Gauss did not change his amiable attitude because she was a woman, eventually, Gauss stopped responding to her letters. Thus, Germain turned back to Legendre for scientific counsel and to share with him her analysis. In 1808, Legendre published the second edition of his {\it Essai sur la th\'eorie des nombres}.

\section{Proving Herself as a Mathematician}

In 1809, Germain shifted her scientific endeavors to mathematical physics, and, in the process, she became the first woman in the history of science to win the coveted grand prix of mathematics awarded by the French Academy of Sciences. Her interest was sparked by the peculiar patterns on vibrating plates, demonstrated by German sound physicist Ernst F. F. Chladni. In 1808, Chladni visited Paris to give talks related to his research in acoustics. His experiments with vibrating surfaces caught the attention of Napoleon, and especially of the scholars at the Academy of Sciences. 

On 27 February 1809, Delambre, secretary of the Academy, announced that {\it Sa Majest\'e Imp\'eriale et Royale} (Napoleon) had approved the idea from the academicians to propose a prize of three thousand francs for the mathematical theory explaining Chladni's experiments. The deadline was set as 1 October 1811.

For the next few years, Germain worked to develop such a theory, reaching out to her mentor, Legendre, as she grappled with Euler's mathematical theories on elastic beams. Her effort culminated with an incomplete differential equation to model a bending plate under static loading. She was the only contestant.  Lagrange, who was one of the reviewers of this work, corrected Germain's results in 1813 by adding the missing term. Lagrange died that same year. After three grueling attempts on her part to match her mathematical theory to Chladni's experiments, while defending ``her biharmonic equation,'' finally in 1816 the Parisian Academy awarded Sophie Germain the grand prize. 

This outcome caused some controversy with Sim\'eon-Denis Poisson, her rival, who tried to diminish her effort. Although not a contestant himself since he was a judge in the competition, Poisson also derived the equation of plate bending on the basis of the theory of elasticity. Other well-known French mathematicians also worked in the same area. After the award, Poisson, Augustin-Louis Cauchy, and Claude-Louis Navier continued expanding this theory, shunning Germain despite her attempts at working with them. Regardless of the controversies and debates that ensued with the top scientists of her time, Sophie Germain was the first person to conceive the theory that would lead to the general plate equation known today as the Germain-Lagrange equation.

This article is about Germain's work in pure mathematics. Hence, the interested reader may find in {\it Sophie Germain}{\footnote{Musielak, D., {\it Sophie Germain, chapters 6, 7}, {\it see} [MD]. }} many details of the intriguing history of her ideas regarding the effect of elastic force on plate vibration.

\section{Returning to Her First Love: Final Attempt to Prove FLT}

On the same day Sophie Germain won the prize for her work on vibrating plates, the commissioners of the Academy of Sciences proposed as a new topic the (last) theorem of Fermat for the new contest of mathematics.{\footnote{{\it S\'eance du mardi 26 d\'ecembre} 1815, pp. 595-596, {\it see} [IP], Tome V.}}  In an effort to bring more honor to France, and at the same time wanting to provide mathematicians with the opportunity to develop ``this part of Science,'' the French academicians issued for the prize of mathematics the general proof of Fermat's Last Theorem. The prize was set as a gold medal valued at three thousand francs, which would be awarded in January 1818.

After the official announcement of the prize, on 26 February 1816, Legendre presented to the Parisian Academy a copy of his newly published Supplement to the second 1808 edition of his {\it Essai}.{\footnote{{\it S\'eance du lundi} 26 f\'evrier 1816. {\it M. Legendre pr\'esente à la Classe un exemplaire de son suppl\'ement à l'Essai sur la th\'eorie des nombres, seconde \'edition}, see [IP] Tome VI.}} 

Did Sophie Germain resume her work on proving Fermat's assertion, aiming at the grand prize? And why not? After all, she had just won a grueling mathematics competition, and the gold medal in her hands must have given her much more confidence. Most importantly, she knew the theorem of Fermat, the particular cases that Euler had completed, and she thought she knew a way to tackle the general proof, just as she had said to Gauss twelve years earlier.

Between 1816 and 1818, others like Germain were also trying to do what Euler had not done. Shortly before the deadline for the contest, on 15 December 1817, the Academy received an anonymous entry, identified with the epigraph {\it Labor omnia vincit}.{\footnote{{\it S\'eance du lundi 15 d\'ecembre 1817. L'Acad\'emie reçoit une pi\'ece avec l'\'epigraphe Labor omnia vincit, sur le th\'eor\`eme de Fermat}. [From Virgil's The Georgics, {\it Un travail acharné vient à bout de tout.}] In English, ``work conquers all'']{\it see} [IP], Tome VII.}}  The Academy also received a second memoir on 26 January 1818, which seems to me it would have been disqualified as it arrived after the deadline. This entry could have been from Sophie Germain. The review panel led by Legendre also included Pierre-Simon de Laplace, Poisson, Louis Poinsot and Cauchy. In any case, on the 9th of March, Legendre reported on the submissions received on the theorem of Fermat, stating that ``These entries contain nothing that deserves a reward.'' Once again, at the public meeting on 16 March 1818, the Academy announced the prize of mathematics to demonstrate ``the remaining theorem of Fermat.'' The deadline for submitting the proof was set as 1 January 1820.

A year later, on 20 February 1819, the Academy received two memoirs sent by Prosper Coste, an engineer graduated from the \'Ecole Polytechnique, one of which addressed, unsuccessfully, the generalization of Fermat's theorem.{\footnote{{\it G\'en\'eralisation du th\'eor\`eme de Fermat sur les doubles \'egalit\'es}, Proc\`es-verbaux des s\'eances de l'Acad\'emie tenues depuis la fondation de l'Institut jusqu'au mois d'ao\^ut 1835. TOME VI, AN 1816-1819, {\it see} [IP].}}  

At that same time, Sophie Germain was busily working on her own proof.{\footnote{Del Centina believes that Germain's work on FLT is datable between 1819 and 1820, {\it see} [DC], p. 350. However, other scholars disagree; for example, {\it see} [L\&P] and their comprehensive assessment of Germain's work.}}

\section{Reaching Out to Gauss, Again}

On 10 May 1819, German-Danish astronomer Heinrich Christian Schumacher wrote to Gauss from Paris to say he would be ``dining with Mademoiselle Germain.'' The visitor at her door must have brought much joy to Sophie Germain. Following their dinner, she entrusted a letter to Schumacher to deliver to Gauss.This is perhaps one of the most significant letters that Germain wrote about her work, containing her strategy for obtaining a general proof of Fermat's Last Theorem.  

However, the scientific contents of the letter remained hidden for over 186 years, until it was carefully evaluated in 2005 by Italian historian A. Del Centina. We don't know what Gauss thought of her approach, but in Section 8 below we shall see the essence of what Sophie Germain wrote to him and her analysis.{\footnote{Del Centina published the letter in its entirety (in its original French) and added clarifying comments, {\it see} [DC], pp. 356-362.}}

At the end of 1819, the Academy of Sciences received more anonymous memoirs to compete for the contest to prove FLT. One memoir identified by the motto {\it Omnia gravia sunt dum ignores; ubi cognoveris, facilia},{\footnote{From Lactantius: ``Those aspiring to perfection, must arrive there through the most unpleasant difficulties.''}} was received on November 22. On 13 December, there were two competing memoirs, and by 27 December one more was submitted with the epigraph {\it Aggrediar non tam perficiendi etc.}.  Laplace, Poisson, Legendre, Poinsot, and Lacroix were to judge the entries. On 24 January 1820, the Commission received an anonymous request, asking for an extension of the term set for the competition. It was denied.

Evidently, none of the competing entries could prove the remaining theorem of Fermat. On 6 March 1820, by recommendation of the Commission of judges, the Academy decided that the prize would not be given and that the topic would be withdrawn from the competition. The sum intended for the prize on Fermat's theorem was given to works in astronomy.{\footnote{The Institute split the 3000 francs between French astronomer Marie-Charles Damoiseau, and Italian astronomers Francesco Carlini and Giovanni Plana, who were awarded the prize for separate memoirs giving new Lunar Tables (positions of the Moon), see [IP], Tome VII.}} 

Five years after the contest to prove FLT ended without a winner, it became known that Sophie Germain had made an important contribution towards the proof. She made no boisterous announcement. There was no public recognition by the academicians either. Rather, Legendre included part of her work in a footnote of his Second Supplement to his {\it Essai} published in September 1825. In this 40-page supplement{\footnote{{\it Essai sur la th\'eorie des nombres, Second suppl\'ement} (Paris, Courcier), septembre 1825, 40 p. (Bibl. de l'Institut, in-4\o, M 572 A(3)).}}, Legendre focused on topics of indeterminate analysis and particularly on Fermat's theorem,{\footnote{In 1827, Legendre published the same work as {\it Recherches sur quelques objets d'analyse ind\'etermin\'ee et particuli\`erement sur le th\'eor\`eme de Fermat}, in the M\'emoires de L'Acad\'emie Royale des Sciences de l'Institut de France, Ann\'ee 1823 [should be 1825].}} where he included the proof of FLT for $n = 5$, with Dirichlet's proof of the same case with an argument for case 2 (an extension of his approach for case 1). 

\section{Sophie Germain Theorem}

After Euler, all that remained was to prove FLT for odd primes $n$. Consider this: if an odd prime $p$ divides $n$, then we can make the reduction $(x^q)^p+(y^q)^p=(z^q)^p$, and redefining the arguments we obtain the Fermat equation (2), $x^p+y^p=z^p$.

If no odd prime divides $n$, then $n$ is a power of 2, so $4|n$ and, in this case, the above equation works with 4 in place of $p$. Since the case $n=4$ was proved by Fermat to have no solutions, it is sufficient to prove the theorem by considering odd prime powers only. Fermat's Last Theorem is traditionally split into two cases:

           Case 1: where exponent $p$ does not divide $xyz$;
           
           Case 2: where exponent $p$ does divide $xyz$.

For odd prime exponent $p$, Case 1 asserts that there do not exist nonzero, pairwise relatively prime integers $x$, $y$ and $z$ such that $x^p+y^p+z^p=0$ and $p$ does not divide $xyz$.

Through a very ingenious analysis, Sophie Germain proved that, if $p$ is an odd prime satisfying two special conditions, then Case 1 of Fermat's theorem holds for $p$. This Germain verified for $p<100$.
 
In the footnote of the art. 22 addressing the two conditions that must be satisfied, Legendre proudly noted that such proposition was due to Germain;{\footnote{{\it Cette d\'emonstration qu'on trouverons tr\`es ing\'enieuse, est due à Mlle Sophie Germain, qui cultive avec succ\`es les sciences physiques et math\'ematiques, comme le prouve le prix qu'elle a remport\'e \`a l'Acad\'emie sur les vibrations des lames \'elastiques. On lui doit encore la proposition de l'art. 13 et celle qui concerne la forme particuli\`ere des diviseurs premiers de  $\alpha$, donnée dans l'art. 11}. [This ingenious demonstration is due to Miss Sophie Germain, who cultivates with success the physical and mathematical sciences, as she proved with the award from the Academy that she won for her work on the vibrations of elastic membranes.  One owes to her also the proposition of article 13 and that one which concerns the particular form of prime divisors of $\alpha$, given in article 11.]. {\it Essai}, Second Supplement, pp. 13-14}} this proposition would later be named Sophie Germain's Theorem. 

Let's now examine Sophie Germain's theorem, using modern language. \\

{\bf Sophie Germain's Theorem}. {\it For an odd prime exponent p, if there exists an auxiliary prime $\theta$ such that there are no two nonzero consecutive $p$th powers modulo $\theta$, nor is $p$ itself a $p$th power modulo $\theta$, then in any solution to the Fermat equation $z^p=x^p+y^p$, one of $x$, $y$, or $z$ must be divisible by $p^2$}.\\

By producing a valid auxiliary prime,{\footnote{Auxiliary primes are prime numbers whose existence is needed to meet the conditions of Germain's proof.}}  Sophie Germain's Theorem can be applied for many prime exponents to eliminate the existence of solutions to the Fermat equation involving numbers not divisible by the exponent $p$. This elimination is today called Case 1 of Fermat's Last Theorem. For example, to prove that any solution to the equation $z^3=x^3+y^3$ (prime exponent $p=3$) would have to have one of $x$, $y$, or $z$ divisible by $p^2=9$, it requires to show that (with auxiliary prime $\theta =13$) no pair of the nonzero cubic residues 1, 5, 8, 12 modulo 13 are consecutive, and $p=3$ is not itself among the residues. 

In the same fashion, we can look for solutions for any other exponent of Fermat's equation. And indeed, this is what Legendre did in the Second Supplement. He verified Germain's hypothesis and generated a table to show, for $p < 100$, the $p$-th power residues modulo $\theta$ for a single auxiliary prime $\theta$ chosen for each value of $p$. He continued to develop more theoretical means of verifying the hypotheses of Germain's Theorem and to demonstrate that any solutions to the Fermat equation for certain exponents would have to be extremely large.{\footnote{For a proof of Germain's theorem, see for example Chapter IV in Ribenboim's {\it Fermat's Last Theorem for Amateurs}, Springer; (1999), or Edwards, {\it Fermat's Last Theorem}, Springer (1997), {\it see} [E].}} 

Sophie Germain's proposition is the first general result about arbitrary exponents for FLT. Her theorem automatically proves Case 1 whenever $2p+1$ is prime. Today, these numbers are known as Germain primes. By definition, a prime number $p$ is a Germain prime if it can be doubled, adding 1 to get another prime, in other words, if $2p+1$ is also prime.

There are many primes $p$ for which $2p+1$ is also prime, but it is unknown whether there are infinitely many such primes. The smallest Germain prime is $p=2$ since 2$\times$2+1=5, which is prime. The next is $p=3$ because 2$\times$3+1=7, also prime, and the next are 5, 11 and 23. Note that Germain primes include 2, the only even prime. The largest known (proven) Sophie Germain prime pair is given by $(p,2p+1)$ where $p=2618163402417$ $\times$ $2^{1290000}-1$, each of which has 388,342 digits.{\footnote{Caldwell, C. K., Prime Pages. The Top Twenty: Sophie Germain. http://primes.utm.edu/top20/page.php?id=2\#records.}}

\section{Sophie Germain's Unexpected Revelation}

Around 1819, Sophie Germain was attempting to prove the impossibility of satisfying Fermat's equation with integers, but she never published her work. Among her papers archived at the Biblioth\`eque Nationale de France, there is a twenty-page manuscript entitled {\it Remarque sur l'impossibilit\'e de satisfaire en nombres entiers a l'\'equation} $x^p+y^p=z^p$. One could speculate that this was among the anonymous memoirs received by the Academy between 1819 and 1820. In 2005, Andrea Del Centina announced the existence of a second copy of this manuscript that he found at the Moreniana Library of Florence (Italy) {\footnote{Del Centina, A., ``Letters of Sophie Germain preserved in Florence.'' His. Math. (2005), 32, pp. 60 - 75.}} and in 2008, Del Centina published a transcript of the manuscript and from it and Germain's last letter to Gauss in 1819 and an appendix she included in the first letter, Del Centina assessed Germain's work on FLT based on this analysis{\footnote{Del Centina, A., ``Unpublished manuscripts of Sophie Germain and a revaluation of her work on Fermat's Last Theorem,'' Arch. Hist. Exact Sci. (2008) 62, pp. 349-392, {\it see} [DC].}} 

In 2010, Reinhard Laubenbacher and David Pengelley [L\&P]{\footnote{Laubenbacher, R. and Pengelley, D.,``{\it Voici ce que j'ai trouv\'e}: Sophie Germain's grand plan to prove Fermat's Last Theorem,'' His. Math. (2010), 37, pp. 641-692, {\it see} [L\&P].}} made a full assessment of Sophie Germain's grand plan as she outlined it.

On 12 May 1819, Sophie Germain wrote to Gauss that this idea [to prove the theorem of Fermat] had been in her mind ever since reading his {\it Disquisitiones}, and immediately she shared the analysis she had carried out. Germain outlined her grand plan for a general proof, including auxiliary theorems and intermediate results that were needed to advance it. Her plan included a unique approach for producing, for each odd prime exponent $p$, an infinite sequence of qualifying auxiliary primes, which would prove the famous theorem.{\footnote{Ibid. pp. 659-661.}} 

Germain wrote that she had recognized a connection between the theory of residues and Fermat's equation, which in fact she had told Gauss about in her correspondence ten years earlier. This letter is dated 12 May 1819, several months before the Academy of Sciences would withdraw the contest for proving Fermat's theorem. Sophie Germain had been working towards that goal and needed a referee, who better than Gauss to assess her proof? Germain asked him for a critique: ``{\it Je vous aurez la plus grande obligation si vous \'et\'es assez bon pour prendre la peine de me dire ce que vous pensez de la marche que j'ai suivie}.''{\footnote{``Unpublished manuscripts,'' p. 362, {\it see} [DC].}}  No letter has been found to determine whether or not Gauss replied.

Independently, Del Centina [DC]{\footnote{Ibid. pp. 356-362.}} and Laubenbacher and Pengelley [L\&P]{\footnote{Laubenbacher and Pengelley, ``{\it Voici ce que j'ai trouv\'e},'' pp. 641-692, {\it see} [L\&P].}} studied this important letter along with others of her undated manuscripts. They reveal to us the depth and scope of the work Sophie Germain carried out in her effort to prove FLT. 

In July 1819, Sophie wrote to Louis Poinsot{\footnote{Del Centina, ``Letters of Sophie,'' p. 62. Others believed that this letter (dated 2 July 2 1819) was addressed to Libri.}}  in reference to work he did in number theory and that she found relevant for her attempted proof of FLT. Poinsot was a professor at the \'Ecole Polytechnique who studied Diophantine equations and found ways to express primes as the difference of two squares and primitive roots. In a paper{\footnote{Poinsot, L., {\it Sur l'alg\`ebre et sur la th\'eorie des nombres}. Read at the Academy of Sciences on 5 May 1817. M\'emoires de la Classe des Sciences Math\'ematiques et Physiques de l'Institut de France, Ann\'ees 1813, 1814, 1815. p. 381, {\it see} [IP] Tome V.}}  presented to the Parisian Academy in 1817, Poinsot developed simultaneously the congruence $x^n \equiv 1$(mod $p$) and the equation $x^n-1=0$ to give an analytic representation of power residues by means of imaginary roots of unity.{\footnote{``Unpublished manuscripts,'' p. 361, {\it see} [DC].}}  Sophie believed that Poinsot's approach was very useful when one is looking for the $p$th power primitive residues.

For Sophie Germain, this epoch was characterized by intense mathematical research and a closer relationship with Legendre. A letter from him dated 31 December 1819 responds to analysis she sent him, addressing an equation which seems to lead nowhere. For some reason, Legendre discouraged Sophie from pursuing this line of research, warning her `` ... {\it que depuis que je vous ai parl\'e pour la premi\`ere fois de ce moyen de recherche, l'opinion que j'avais qu'il pouvait r\'eussir est maintenant bien affaiblie, et qu'en somme je crois qu'il sera aussi st\'erile que bien d'autres}.''{\footnote{Stupuy, {\it Ouvres Philosophiques}, pp. 310-311, {\it see} [SH].}}  Sophie was not easily discouraged.	\\

{\bf Germain's Research to Prove Fermat's Last Theorem}\\

When Legendre published a third edition of his {\it Essai}, renamed {\it Th\'eorie des nombres}, he did not mention Sophie Germain's extensive work to prove FLT. This omission helps us to understand why, for almost two centuries, it was assumed that her contribution to this famous proof was just as Legendre summarized it in the footnote of his Second Supplement in 1825. However, the reevaluation of her manuscripts and her correspondence with Legendre and Gauss by DC and L\&P indicates otherwise. They carried out a comprehensive study of Sophie Germain's work, publishing in 2010 a reassessment of her ``grand plan'' for proving Fermat's Last Theorem. 

The scholarly work by Laubenbacher and Pengelley provides the most comprehensive evaluation of Sophie Germain's contribution to proving FLT to date. They discovered a wealth of results in these manuscripts beyond the single theorem for Case 1 for which Germain was known. They analyzed the supporting algorithms that she invented and found that they were based on ideas and results that other number theorists discovered independently much later.{\footnote{Laubenbacher and Pengelley, ``{\it Voici ce que j'ai trouv\'e},'' {\it see} [L\&P].}}  

Sophie Germain's Theorem is simply a small part of her big program, a piece that could be applied separately as an independent theorem. Germain's objective was to prove FLT for exponent $p$ by producing an infinite sequence of qualifying auxiliary primes of the form $\theta =2Np+1$. She developed methods to validate a qualifying requirement for infinitely many auxiliary primes, called the Non-Consecutivity (N-C) condition: {\it There do not exist two nonzero consecutive $p$-th power residues, modulo $\theta$}.

Starting from Fermat's equation $x^p+y^p=z^p$, $p$ an odd prime number, Germain claimed that if this equation were possible, then every prime number of the form $2Np+1$ ($N$ being any integer), for which there are no two consecutive $p$-th power residues in the sequence of natural numbers, necessarily divides one of the numbers $x$, $y$, and $z$. To prove FLT, Germain intended to prove that if there were infinitely many such qualifying numbers $2Np+1$, then the Fermat equation would be impossible. So Germain's plan for proving FLT for exponent $p$ required that she develop methods to validate the qualifying condition for infinitely many auxiliary primes of the form $\theta =2Np+1$. Sophie Germain developed an algorithm verifying the condition within certain ranges, and outlined an induction on auxiliaries to carry her plan forward.{\footnote{Ibid. pp. 657, 682.}} 

To establish this condition for various $N$ and $p$, she carried out her analysis over many pages, including the general consequences of nonzero consecutive $p$th power residues modulo a prime $\theta=2Np+1$, where $N$ is never a multiple of 3. Germain included a summary table with her results verifying Condition N-C for auxiliary primes $\theta$ using relevant values of $N \leq 10$ and primes $2<p<100$, confident that she could easily extend the range.

Only a number theorist can grasp the sophisticated plan that Sophie Germain conceived to prove Fermat's Last Theorem, which, as noted by L\&P,  used the multiplicative structure in a cyclic prime field and a set (group) of transformations of consecutive $p$th powers. However, Germain's plan cannot be carried to completion, as we will discuss shortly. Moreover, L\&P also found an undated letter to Legendre in which Germain actually proved that her plan fails for $p=3$: For any prime $\theta$ of the form $6a+1$, with $\theta >13$, there are (nonzero) consecutive cubic residues. In other words, the N-C condition fails for $\theta =2Np+1$ when $p=3$ and $N>2$, so the only valid auxiliary primes for $p = 3$ for the N-C condition are $\theta =7$ and 13.

Sophie Germain also attempted to prove FLT by related means for other special forms of the exponent, but none of these approaches succeeded.{\footnote{Ibid. p. 679.}}  She knew that she had not proved the theorem of Fermat even for a single exponent. Hence, she wanted to show for specific exponents that any possible solutions to the Fermat equation would have to be extremely large. In her letter to Gauss, she stated that any possible solutions would consist of numbers ``whose size frightens the imagination.'' Germain stated, proved, and applied her ``large size theorem'' to show for $p=5$, the valid auxiliary primes $\theta =11,41,71,101$ indicate that any solution to the Fermat equation would force a solution number to have at least 39 decimal digits. Gauss could have considered her argument and perhaps seen that this theorem and its applications were not valid. But it appears that he did not reply to her letter. 

The details of Sophie Germain's analysis and a comprehensive evaluation of her proofs are given by L\&P. They elucidated her argument, translated her analysis, and verified her claim.{\footnote{Ibid. p. 664.}}  Laubenbacher and Pengelley found her proof rather remarkable, especially since she performed it over the course of one night. As they concluded, Sophie Germain was an impressive number theorist.

Researchers have asked, ``Did Germain, in the end, become aware of the futility of her plan?'' I believe that she did. Sophie Germain knew that her plan for proving FLT was hopeless. But when did she find this out? In his 1829 memoir, Guglielmo Libri asserted that seeking infinitely many auxiliary primes would not work. He proved for $n=3$ and $n=4$ and stated for any $n$ that for the primes $p$, exceeding an assignable limit, Germain's congruence has solutions prime to $p$, so that it is futile to attempt to prove Fermat's theorem by trying to show that one of the unknowns is divisible by an infinitude of primes. Libri claimed that this result was already in two of his earlier memoirs presented to the Parisian Academy of Sciences in 1823 and 1825. From this claim, it was concluded that close followers of the Academy should have been aware by 1825 that Libri's work would doom the auxiliary prime approach to FLT.{\footnote{Ibid. p. 662.}}  But, is this really how it happened?

There is no record of a review of Libri's work in 1823. On 9 August 1824, during a meeting of the Academy, Cauchy gave a verbal report on Libri's memoir,{\footnote{M\'emoires de l'Acad\'emie des sciences de l'Institut de France, Tome VII, Ann\'ee 1824. Chez Firmin Didot, p\`ere et fils, libraries, Paris, 1827. pp. lxviii-lxix, {\it see} [MIF].}}  which was divided into three parts. In part 1, Libri showed a method to provide whole solutions of indeterminate equations. In part 2, titled ``theory of congruences,'' Libri determined the number of solutions of a congruence on the basis of these same roots. In part 3, Libri established analytical formulas with which he expressed the number of solutions of an indeterminate equation.{\footnote{Tome VIII, Années 1824-1827, pp. 121-123, {\it see} [IP].}}  

From Cauchy's overview we can infer that Libri established ``in a very simple way Euler's formula relating to the divisors of numbers.'' Libri also proved that the relationships between the coefficients of algebraic equations and their roots extend to congruences in which all roots are real. He inferred from this principle Wilson's theorem, which states that a natural number $n > 1$ is a prime number if and only if $(n-1)! \equiv -1$ (mod $n$). In other words, Wilson's theorem asserts that the factorial is one less than a multiple of $n$ exactly when $n$ is a prime number.

Then again, Cauchy did not mention FLT or its proof, nor did he note that Libri had determined that ``for exponents 3 and 4, there can be at most finitely many auxiliary primes satisfying the non-consecutive (N-C) condition.'' Cauchy would have emphasized such a crucial result if indeed it were in that memoir. Even if he had obtained it, as Libri claimed, Sophie Germain would not have known about that result before 1824. She met Libri in May of 1825.

During his first visit to Paris, Libri presented his work to the Academy at a meeting{\footnote{Ibid., p. 223.}} on 13 June 1825.  Again, he made no reference to analysis related to auxiliary primes. A month later, German mathematician Gustav Lejeune Dirichlet submitted to the Paris Academy his famous {\it M\'emoire sur l'Impossibilit\'e de quelques \'equations ind\'etermin\'ees du 5e degr\'e},{\footnote{Ibid., p. 239.}} in which he tackled FLT for exponent $n=5$. 

The proof of Fermat's theorem was an exciting topic of discussion among academicians when Libri and Sophie Germain met in May 1825. It is conceivable that during their scientific exchange, Germain shared with him the results of her attempts at the proof. The existence of some of Sophie Germain's manuscripts in an Italian library today suggests that Libri had them in his possession. After he left France, Libri wrote to Germain asking her to speak to Joseph Fourier in order to expedite the report on the memoir that he had submitted six months earlier.{\footnote{Del Centina, ``Letters of Sophie Germain preserved in Florence.'' Hist. Math. 32 (2005).}} 	

On 13 March 1826, Fourier and Cauchy made a verbal report of Libri's paper.{\footnote{Proc\`es-verbaux des s\'eances de l'Acad\'emie des sciences, Tome VIII, Ann\'ees 1824-1827, p. 358, {\it see} [IP].}}  Cauchy stated that Libri ``had shown the method by which he determined equal roots of numerical equations that contain a single unknown and can be extended to an indeterminate equation; and then he established principles specifically applicable to equations that are linear over one of the unknowns, and that Gauss named congruences.''{\footnote{M\'emoires de l’Acad\'emie des sciences de l'Institut de France, Tome IX, Ann\'ee 1826. Chez Firmin Didot Fr\`eres, libraries, Paris, 1830.}}  Again, this report does not assert that the attempts of others to prove FLT by finding infinitely many such auxiliaries are in vain. This conclusion would have made a huge impact in the work of Legendre and Germain, and Legendre would have noted it in his last paper{\footnote{Legendre, M.-A., {\it M\'emoire sur la d\'etermination des fonctions Y et Z qui satisfont \`a l'\'equation  $4(Xn-1)=(X-1)(Y^2 \pm nZ^2)$, $n$ \'etant un nombre premier $4i \mp 1$}. M\'emoires de l'Acad\'emie des Sciences, Vol. XI, Année 1828, pp. 81-99.}} published in 1828.

I find it rather strange that when Legendre published the third edition of his {\it Th\'eorie des nombres}, he did not include the method that Sophie Germain had pioneered.{\footnote{Laubenbacher and Pengelley, ``{\it Voici ce que j'ai trouv\'e},'' pp. 641-92, {\it see} [L\&P].}}  More intriguing perhaps is why Germain did not publish her theorem or her extensive results to prove FLT. Maybe she was waiting for validation, especially from Gauss, the scholar she had asked repeatedly for a critique of her work and who did not tell her what he thought of her proof. Or perhaps, she knew that her effort was futile and that the general proof was out of reach. Nevertheless, the elegant theorem that Sophie Germain developed is an important contribution to the proof of Fermat's Last Theorem.

\section{Epilogue}

Sophie Germain stood in the midst of the so-called Golden Age of number theory, a period of time linking Euler to Kummer.{\footnote{Edwards, {\it Fermat's Last Theorem.}, p. 61, {\it see} [E].}} Attempting to prove Fermat's Last Theorem gave inspiration to many mathematicians and led them to develop new ideas. These mathematical ideas include Ernst Kummer's theory of ideal factors, a concept he developed in order to deal with the higher reciprocity laws from which the general theory of algebraic numbers grew. It also includes Dirichlet's analytical formula, which is used for the number of classes of binary quadratic forms with given determinant.{\footnote{Ibid. p. 60.}}  

Fermat's Last Theorem became one of the most alluring mathematical assertions in the history of mathematics.  By the twentieth century, Fermat's assertion had been verified for exponents up into the thousands, and this increased when the calculations were facilitated by the use of computers. However, not one person could prove that FLT is true for an infinite set of prime exponents. Clearly, rather than verifying (1) case by case up to $n$ approaching infinity, a new approach had to be developed in order to obtain an absolute proof of FLT for all $n>2$. 

Finally, in 1995, British mathematician Andrew Wiles did just that. His proof combined results from several concepts of modern mathematics, relying heavily on another conjecture, an idea developed in the 1950’s that connects topology and number theory. Japanese mathematicians Taniyama and Shimura supposed that ``every elliptic curve over the rationals is modular.'' If true, their conjecture would imply FLT. But the Taniyama-Shimura conjecture, now known as the modularity theorem,{\footnote{http://mathworld.wolfram.com/Taniyama-ShimuraConjecture.html }} was just as challenging to prove as FLT itself.

At a conference in Cambridge in 1993, Wiles{\footnote{At that time, Andrew Wiles was a professor at Princeton University.}} gave a series of lectures on ``Modular forms, elliptic curves and Galois representations.'' After presenting many new results, Wiles concluded by stating: ``This proves Fermat's Last Theorem.'' The awe and excitement among those present must have been palpable, and the news of such a mathematical discovery quickly spread and appeared on the front pages of newspapers across the world. Of course, as any mathematical proof, the new results had to be subjected to an intense scrutiny by specialized mathematicians. 

The referees soon discovered a flaw in the proof. To resolve it, Wiles asked help from another English mathematician, Richard Taylor. After months of intense work, both Wiles and Taylor succeeded. Their results appeared in two papers in the 1995 Annals of Mathematics, volume 141. The main paper is by Wiles (109 pages long);{\footnote{Wiles, A., {\it Modular elliptic curves and Fermat's Last Theorem}. Annals of Mathematics (1995) 141 (3): 443–551. doi:10.2307/2118559. JSTOR 2118559. OCLC 37032255.}}  its companion is by Taylor and Wiles (20 pages long){\footnote{Taylor, R. and Wiles, A., {\it Ring theoretic properties of certain Hecke algebras}. Annals of Mathematics (1995) 141 (3): 553–572. doi:10.2307/2118560. JSTOR 2118560. OCLC 37032255.}}  to establish results that are essential for the proof.

Considering the number of pages required for Wiles and Taylor to prove FLT, perhaps this is what Fermat meant when he scribbled that ``This margin is too narrow to contain [his] marvelous proof.'' It is evident that Fermat's method would have been much different from the approach taken by Wiles.

The mystery that shrouds the work of Sophie Germain remains. Did she do more than her manuscripts reveal? She came to realize that her original plan for proving FLT was hopeless. I firmly believe that Sophie Germain reached the conclusion that her analysis lacked the necessary tools, and that this deficiency trampled her effort to prove {\it la c\'el\'ebre \'equation de Fermat}. We may never know for sure ..., unless someone finds more of her personal manuscripts hidden in a library vault somewhere. 

Sophie Germain died in Paris at one in the morning, 27 June 1831.\\

{\bf Acknowledgement}\\

I am deeply indebted to Professor David Pengelley for reading the entire manuscript and providing me with much needed feedback. He challenged some of my original statements, forcing me to reevaluate my sources; his critique ensured that every aspect of Sophie Germain's contribution to mathematics that I emphasize in this paper is in fact fully validated by historical evidence and scholarly assessments. Any remaining faults are entirely my own, of course. 

Thank you, David. {\it Quod postquam tandem ex voto successisset, illecebris harum quaestionum ita fui implicatus, ut eas deserere non potuerim.}\\

{\bf References}\\

[DC] Del Centina, A. ``Unpublished manuscripts of Sophie Germain and a revaluation of her work on Fermat's Last Theorem,'' Arch. Hist. Exact Sci. (2008) 62.

[E] Edwards, H. M., {\it Fermat's Last Theorem. A Genetic Introduction to Algebraic Number} Theory, Springer-Verlag, New York (1997)

[F] Fuss, P.-H., ed., {\it Correspondance Math\'ematique et Physique}, Imp. Acad. Sci., St. Petersburg, 1843, vol. I. 

[IP] Institut de France, Acad\'emie des Sciences. Proc\`es-verbaux des S\'eances de l'Acad\'emie tenues depuis la fondation de l'Institut. Tome I, An IV-VII (1795-1799); Tome IV, An 1808-1811; Tome V, An 1812-1815; Tome VI, An 1816-1819; Tome VII, An 1820-1823; Tome VIII, An 1824-1827.

[L\&P] Laubenbacher, R. and Pengelley, D., 2010. ``{\it Voici ce que j'ai trouv\'e}: Sophie Germain's grand plan to prove Fermat's Last Theorem,'' Historia Mathematica, Vol. 37, pp. 641-692.

[MD] Musielak, D., {\it Sophie Germain: Revolutionary Mathematician}. Springer Biography, ISBN 978-3-030-38374-9, (2020) 

www.springer.com/gp/book/9783030383749

[MIF] {\it M\'emoires de l'Acad\'emie des sciences de l'Institut de France, Ann\'ee 1823, Tome VI, }. Chez Firmin Didot Fr\`eres, libraries, Paris (1827)

[SH] Stupuy, H., {\it Ouvres philosophiques de Sophie Germain, suivies de pens\'ees et de lettres in\'edites et pr\'ec\'ed\'ees d'une notice sur sa vie et ses ouvres}, Paris, P. Ritti (1879) in-18.\\

$\hskip0.7in$                              {Dora E Musielak; University of Texas at Arlington, June 2019}

\end{document}